# Learning gradients on manifolds

SAYAN MUKHERJEE[1], QIANG WU[2] and DING-XUAN ZHOU[3]

[1]*Department of Statistical Science and Department of Mathematics Institute for Genome Sciences & Policy, Department of Computer Science, Department of Biostatistics and Bioinformatics, Duke University, Durham, NC 27708, USA. E-mail: sayan@stat.duke.edu*
[2]*Department of Mathematics, Michigan State University, East Lansing, MI 48824, USA. E-mail: wuqiang@math.msu.edu*
[3]*Department of Mathematics, City University of Hong Kong, Tat Chee Avenue, Kowloon, Hong Kong, China. E-mail: mazhou@cityu.edu.hk*

A common belief in high-dimensional data analysis is that data are concentrated on a low-dimensional manifold. This motivates simultaneous dimension reduction and regression on manifolds. We provide an algorithm for learning gradients on manifolds for dimension reduction for high-dimensional data with few observations. We obtain generalization error bounds for the gradient estimates and show that the convergence rate depends on the intrinsic dimension of the manifold and not on the dimension of the ambient space. We illustrate the efficacy of this approach empirically on simulated and real data and compare the method to other dimension reduction procedures.

## 1. Introduction

The inference problems associated with high-dimensional data offer fundamental challenges – the scientifically central questions of model and variable selection – that lie at the heart of modern statistics and machine learning. A promising paradigm in addressing these challenges is the observation or belief that high-dimensional data arising from physical or biological systems can be effectively modeled or analyzed as being concentrated on a low-dimensional manifold. In this paper we consider the problem of dimension reduction – finding linear combinations of salient variables and estimating how they covary – based upon the manifold paradigm. We are particularly interested in the high-dimensional data setting, where the number of variables is much greater than the number of observations, sometimes called the "large $p$, small $n$" paradigm [22].

The idea of reducing high-dimensional data to a few relevant dimensions has been extensively explored in statistics, computer science and various natural and social sciences. In machine learning the ideas in isometric feature mapping (ISOMAP) [20], local linear embedding (LLE) [18], Hessian eigenmaps [8] and Laplacian eigenmaps [2] are all formulated from the manifold paradigm. However, these approaches do not use response variates in the models or algorithms and hence may be suboptimal with respect to predicting response. In statistics the ideas developed in sliced inverse regression (SIR) [13], (conditional) minimum average variance estimation (MAVE) [23] and sliced average variance estimation (SAVE) [6] consider dimension reduction for regression problems. The response variates are taken into account and the focus is on predictive linear subspaces called either effective dimension reduction (e.d.r.) space [23] or central mean subspace [5]. These approaches do not extend to the manifold paradigm. In [15, 16] the method of learning gradients was introduced for variable selection and regression in high-dimensional analysis for regression and binary regression setting. This method, in machine learning terminology, is in spirit a supervised version of Hessian eigenmaps and, in statistics terminology, can be regarded as a nonparametric extension of MAVE. This approach can be extended to the general manifold setting. The main purpose of this paper is to explore this idea.

The inference problem in regression is estimating the functional dependence between a response variable $Y$ and a vector of explanatory variables $X$

$$Y = f(X) + \epsilon$$

from a set of observations $D = \{(x_i, y_i)\}_{i=1}^n$ where $X \in \mathcal{X} \subset \mathbb{R}^p$ has dimension $p$ and $Y \in \mathbb{R}$ is a real valued response for regression and $Y \in \{\pm 1\}$ for binary regression. Typically the data are drawn i.i.d. from a joint distribution, $(x_i, y_i) \sim \rho(X, Y)$. We may in addition want to know which variables of $X$ are most relevant in making this prediction. This can be achieved via a variety of methods [4, 12, 21]. Unfortunately, these methods and most others do not provide estimates of covariance for salient explanatory variables and cannot provide the e.d.r. space or central mean subspace. Approaches such as SIR [13] and MAVE [23] address this shortcoming.

SIR and its generalized versions have been successful in a variety of dimension reduction applications and provide almost perfect estimates of the e.d.r. spaces once the design conditions are fulfilled. However, the design conditions are limited and the method fails when the model assumptions are violated. For example, quadratic functions or between group variances near zero violate the model assumptions. In addition, since SIR finds only one direction, its applicability to binary regression is limited.

MAVE and the outer product of gradient (OPG) method [23] are based on estimates of the gradient outer product matrix either implicitly or explicitly. They estimate the central mean subspace under weak design conditions and can capture all predictive directions. However, they cannot be directly used for "large $p$, small $n$" setting due to overfitting. The learning gradient method in [15, 16] estimates the gradient of the target function by nonparametric kernel models. It can also be used to compute the gradient outer product matrix and realize the estimation of the central mean subspace by the same manner as OPG (see Section 4 below). Moreover, this method can be directly used for the "large $p$, small $n$" setting because the regularization technique prevents overfitting and guarantees stability.

All the above methods have been shown to be successful by simulations and applications. However we would like a theoretical and conceptual explanation of why this approach to dimension reduction is successful with very few samples and many dimensions. Conceptually this reduces to the following analysis: For a target function on a nonlinear manifold, the gradient outer product matrix defined in the Euclidean can still be used to estimate predictive directions even when the gradient is not well defined on the ambient space. Theoretically, we notice that the consistency results for MAVE and OPG [23] and learning gradients [15, 16] provide asymptotic rates of order $O(n^{-1/p})$. Clearly this is not satisfactory and does not support practical applicability when $p$ is large, especially for the setting where $p \gg n$. Intuitively one should expect that the rate would be a function not of the dimension of the ambient space but of the intrinsic dimension of the manifold.

In this paper we extend the learning gradient algorithms from the Euclidean setting to the manifold setting to address these questions. Our two main contributions address the conceptual and theoretical issues above. From a conceptual perspective we will design algorithms for learning the gradient along the manifold. The algorithm in the Euclidean setting can be applied without any modifications to the manifold setting. However, the interpretation of the estimator is very different and the solutions contain information on the gradient of the target function along the manifold. This interpretation provides a conceptual basis for using the usual $p$-dimensional gradient outer product matrix for dimension reduction. From a theoretical perspective, we show that the asymptotic rate of convergence of the gradient estimates is of order $O(n^{-1/d})$ with $d$ being the intrinsic dimension of the manifold. This suggests why in practice these methods perform quite well, since $d$ may be much smaller than $n$ though $p \gg n$.

The paper will be arranged as follows. In Section 2 we develop the learning gradient algorithms on manifolds. The asymptotic convergence is discussed in Section 3, where we show that the rate of convergence of the gradient estimate is of order $O(n^{-1/d})$. In Section 4 we explain why dimension reduction via gradient estimates has a solid conceptual basis in the manifold setting and discuss relations and comparisons to existing work. Simulated and real data are used in Section 5 to verify our claims empirically and closing remarks and comments are given in Section 6.

## 2. Learning gradients

In this section we first review the gradient estimation method on Euclidean spaces proposed in [15, 16]. Then after a short discussion of Taylor series expansion on manifolds we formulate learning gradients under the manifold setting.

### 2.1. Learning gradients in Euclidean space

In the standard regression problem the target is the regression function defined by the conditional mean of $Y|X$, that is, $f_r = \mathbb{E}_Y[Y|X]$. The objective of learning gradients is

to estimate the gradient

$$\nabla f_r = \left(\frac{\partial f_r}{\partial x^1}, \ldots, \frac{\partial f_r}{\partial x^p}\right)^{\mathrm{T}},$$

where $x = (x^1, \ldots, x^p) \in \mathbb{R}^p$. The learning gradient algorithm developed in [16] is motivated by fitting first-order differences.

Recall $f_r$ is the minimizer of the variance or mean square error functional,

$$f_r(x) = \mathbb{E}(Y|X=x) = \arg\min \mathrm{Var}(f) \qquad \text{where } \mathrm{Var}(f) = \mathbb{E}(Y - f(X))^2.$$

When only a set of samples $D$ are available the functional is usually approximated empirically

$$\mathrm{Var}(f) \approx \frac{1}{n}\sum_{i=1}^{n}(y_i - f(x_i))^2.$$

Using the first-order Taylor series expansion approximating a smooth function $f$ by

$$f(x) \approx f(u) + \nabla f(x) \cdot (x - u) \qquad \text{for } x \approx u,$$

the variance of $f$ may be approximated as

$$\mathrm{Var}(f) \approx \frac{1}{n^2}\sum_{i,j=1}^{n} w_{ij}[y_i - f(x_j) - \nabla f(x_j) \cdot (x_i - x_j)]^2, \tag{2.1}$$

where $w_{ij}$ is a weight function that ensures the locality of $x_i \approx x_j$ and thus $w_{ij} \to 0$ as $\|x_i - x_j\| \to \infty$. The weight function $w_{ij}$ is typically characterized by a bandwidth parameter, for example, a Gaussian with the bandwidth as the standard deviation $w_{ij} = \mathrm{e}^{-\|x_i - x_j\|^2/(2s^2)}$.

Learning gradient algorithms were specifically designed for very high-dimensional data but with limited number of observations. For regression the algorithm was introduced in [16] by nonparametric reproducing kernel Hilbert space (RKHS) models. The estimate of the gradient is given by minimizing (2.1) with regularization in an RKHS

$$\vec{f}_D := \arg\min_{\vec{f} \in \mathcal{H}_K^p}\left[\sum_{i,j=1}^{n} w_{ij}(y_j - y_i - \vec{f}(x_i) \cdot (x_j - x_i))^2 + \lambda\|\vec{f}\|_K^2\right], \tag{2.2}$$

where $\mathcal{H}_K = \mathcal{H}_K(\mathcal{X})$ is a reproducing kernel Hilbert space (RKHS) on $\mathcal{X}$ associated with a Mercer kernel $K$ (for the definition and properties of RKHS, see [1]) and $\mathcal{H}_K^p$ is the space of $p$ functions $\vec{f} = (f_1, \ldots, f_p)$ where $f_i \in \mathcal{H}_K$, $\|\vec{f}\|_K^2 = \sum_{i=1}^{p}\|f_i\|_K^2$ and $\lambda > 0$.

With the weight function chosen to be the Gaussian with standard variance $s^2$, a finite sample probabilistic bound for the distance between $\vec{f}_D$ and $\nabla f_r$ is provided in [16], which implies the convergence of the gradient estimate to the true gradient, $\vec{f}_D \to \nabla f_r$, at a slow rate, $\mathrm{O}(n^{-1/p})$.

For binary classification problems where $Y = \{\pm 1\}$, the learning gradient algorithm was introduced in [15]. The idea is to use the fact that the function

$$f_c(x) := \log\left[\frac{\text{Prob}(y=1|x)}{\text{Prob}(y=-1|x)}\right] = \log\left[\frac{\rho(y=1|x)}{\rho(y=-1|x)}\right]$$

is given by

$$f_c = \arg\min \mathbb{E}\phi(Yf(X))$$

with $\phi(t) = \log(1 + e^{-t})$. Notice that the Bayes optimal classification function is given by $\text{sgn}(f_c)$, the sign of $f_c$. In the binary classification setting we learn the gradient of $f_c$. Applying the first-order Taylor expansion of $f$ we have for the given data $D$

$$\mathbb{E}\phi(Yf(X)) \approx \frac{1}{n^2} \sum_{i,j=1}^{n} w_{ij}\phi(y_i(f(x_j) + \nabla f(x_i) \cdot (x_i - x_j))).$$

Modeling $f$ and $\nabla f$ by a real valued function $g$ and a vector valued function $\vec{f}$, respectively, leads to the empirical risk

$$\mathcal{E}_{\phi,D}(g, \vec{f}) = \frac{1}{n^2} \sum_{i,j=1}^{n} w_{ij}\phi(y_i(g(x_j) + \vec{f}(x_i) \cdot (x_i - x_j))).$$

Minimizing this empirical risk with a regularization term gives the following algorithm

$$(g_{\phi,D}, \vec{f}_{\phi,D}) = \arg\min_{(g,\vec{f}) \in \mathcal{H}_K^{p+1}} (\mathcal{E}_{\phi,D}(g, \vec{f}) + \lambda_1 \|g\|_K^2 + \lambda_2 \|\vec{f}\|_K^2), \tag{2.3}$$

where $\lambda_1, \lambda_2$ are the regularization parameters. A finite sample probabilistic bound for the distance from $g_{\phi,D}$ to $f_c$ and $\vec{f}_{\phi,D}$ to $\nabla f_c$ is provided in [15], which leads to a very slow rate of order, $O(n^{-1/p})$.

## 2.2. Gradients and Taylor expansions on Riemannian manifolds

In order to extend learning gradients to the manifold setting, it is necessary to formulate gradients and first-order Taylor expansions on manifolds. To do this we need to introduce some concepts and notation from Riemannian geometry. We introduce only what is needed so that we can stress concepts over technical details. For a complete and rigorous formulation, see [7].

The two key concepts are vector fields and the exponential map. Let $\mathcal{M}$ be a $d$-dimensional smooth (i.e., $C^\infty$) Riemannian manifold and $d_\mathcal{M}(a, b)$ be the Riemannian distance on $\mathcal{M}$ between two points $a, b \in \mathcal{M}$. The tangent space at a point $q \in \mathcal{M}$ is a $d$-dimensional linear space and will be denoted by $T_q\mathcal{M}$. There exists an inner product on this tangent space $\langle \cdot, \cdot \rangle_q$ that defines the Riemannian structure on $\mathcal{M}$.

A vector field on a manifold is an assignment to every point $q$ on the manifold tangent vector in $T_q\mathcal{M}$. The gradient of a smooth function $f$ on $\mathcal{M}$, $\nabla_{\mathcal{M}} f$, is a vector field satisfying

$$\langle \nabla_{\mathcal{M}} f(q), v \rangle_q = v(f) \qquad \text{for all } v \in T_q\mathcal{M}, q \in \mathcal{M}.$$

It can be represented using an orthonormal basis $\{e_1^q, \ldots, e_d^q\}$ of $T_q\mathcal{M}$ as the vector

$$\nabla_{\mathcal{M}} f(q) = (e_1^q f, \ldots, e_d^q f).$$

If the manifold is the Euclidean space ($\mathcal{M} = \mathbb{R}^d$), then one may take $e_i^q = \frac{\partial}{\partial q^i}$ and the above definition reduces to the standard definition of gradients.

The exponential map at a point $q$, denoted by $\exp_q$, is a map from the tangent space $T_q\mathcal{M}$ to the manifold $\mathcal{M}$. It is defined by the the locally length-minimizing curve, the so-called geodesic. The image of $v \in T_q\mathcal{M}$ is the end of a geodesic starting at $q$ with velocity $v$ and time 1. In general the exponential map is only locally defined in that it maps a small neighborhood of the origin in $T_q\mathcal{M}$ to a neighborhood of $q$ on the manifold. Its inverse, $\exp_q^{-1}$, maps the point $\exp_q(v)$ to the vector $(v^1, \ldots, v^d) \in \mathbb{R}^d$ where $v = \sum_{i=1}^d v^i e_i^q$. This provides a local chart for the neighborhood of $q$ and $\{e_i^q\}$ are called the $q$-normal coordinates of this neighborhood.

Under the $q$-normal coordinates the gradient vector field $\nabla_{\mathcal{M}} f$ takes the form $\nabla_{\mathcal{M}} f(q) = \nabla(f \circ \exp_q)(0)$. Note that $f \circ \exp_q$ is a smooth function on $T_q\mathcal{M} \cong \mathbb{R}^d$. The following first-order Taylor series expansion holds:

$$(f \circ \exp_q)(v) \approx (f \circ \exp_q)(0) + \nabla(f \circ \exp_q)(0) \cdot v \qquad \text{for } v \approx 0.$$

This gives us the following Taylor expansion of $f$ around a point $q \in \mathcal{M}$:

$$f(\exp_q(v)) \approx f(q) + \langle \nabla_{\mathcal{M}} f(q), v \rangle_q \qquad \text{for } v \in T_q\mathcal{M}, v \approx 0. \tag{2.4}$$

The above expansion does not depend on the choice of the coordinate system at $q$.

### 2.3. Learning gradients on Riemannian manifolds

In the manifold setting, the explanatory variables are assumed to concentrate on an unknown $d$-dimensional Riemannian manifold $\mathcal{M}$ and there exists an isometric embedding $\varphi: \mathcal{M} \to \mathbb{R}^p$ with every point on $\mathcal{M}$ described by a vector in $\mathbb{R}^p$. When a set of points are drawn from the marginal distribution $\rho_{\mathcal{M}}$ concentrated on $\mathcal{M}$ what we know are not the points $\{q_i\}_{i=1}^n \in \mathcal{M}$ themselves but their images under $\varphi: x_i = \varphi(q_i)$.

To formulate the learning gradient algorithm for the regression setting, we apply the first-order Taylor expansion (2.4). The empirical approximation of the variance by the data $\{(q_i, y_i)\}$ is

$$\text{Var}(f) \approx \frac{1}{n^2} \sum_{i,j=1}^n w_{ij} [y_j - y_i - \langle \nabla_{\mathcal{M}} f(q_i), v_{ij} \rangle_{q_i}]^2, \tag{2.5}$$

where $v_{ij} \in T_{q_i}\mathcal{M}$ is the tangent vector such that $q_j = \exp_{q_i}(v_{ij})$. One may immediately notice the difficulty that $v_{ij}$ is not computable without knowing the manifold. A natural idea to overcome this difficulty is to represent all quantities in $\mathbb{R}^p$. This is also compatible with the fact that we are given images of the points $x_i = \varphi(q_i) \in \mathbb{R}^p$ rather than the $d$-dimensional representation on the manifold.

Suppose $x = \varphi(q)$ and $\xi = \varphi(\exp_q(v))$ for $q \in \mathcal{M}$ and $v \in T_q\mathcal{M}$. Since $\varphi$ is an isometric embedding, i.e., $\mathrm{d}\varphi_q : T_q\mathcal{M} \to T_x\mathbb{R}^p \cong \mathbb{R}^p$ is an isometry for every $q \in \mathcal{M}$, the following holds:

$$\langle \nabla_\mathcal{M} f(q), v \rangle_q = \mathrm{d}\varphi_q(\nabla_\mathcal{M} f(q)) \cdot \mathrm{d}\varphi_q(v),$$

where

$$\mathrm{d}\varphi_q(v) \approx \varphi(\exp_q(v)) - \varphi(q) = \xi - x \qquad \text{for } v \approx 0.$$

Applying these relations to the observations $D = \{(x_i, y_i)\}_{i=1}^n$ yields

$$\langle \nabla_\mathcal{M} f(q_i), v_{ij} \rangle_{q_i} \approx \mathrm{d}\varphi_{q_i}(\nabla_\mathcal{M} f(q_i)) \cdot (x_j - x_i). \tag{2.6}$$

Notice further that $\mathrm{d}\varphi \circ \nabla_\mathcal{M} f$ is a $p$-dimensional vector valued function on $\varphi(\mathcal{M}) \subset \mathbb{R}^p$ defined by $(\mathrm{d}\varphi \circ \nabla_\mathcal{M} f)(\varphi(q)) = \mathrm{d}\varphi_q(\nabla_\mathcal{M} f(q))$ for $q \in \mathcal{M}$. Applying (2.6) to (2.5) and denoting $\mathrm{d}\varphi \circ \nabla_\mathcal{M} f$ by $\vec{f}$ yields

$$\mathrm{Var}(f) \approx \mathcal{E}_D(\vec{f}) = \frac{1}{n^2} \sum_{i,j=1}^n w_{ij}[y_j - y_i - \vec{f}(x_i) \cdot (x_j - x_i)]^2.$$

Minimizing this quantity leads to the following learning gradient algorithm on manifolds:

**Definition 2.1.** *Let $\mathcal{M}$ be a Riemannian manifold and $\varphi : \mathcal{M} \to \mathbb{R}^p$ be an isometric embedding that is unknown. Denote $\mathcal{X} = \varphi(\mathcal{M})$ and $\mathcal{H}_K = \mathcal{H}_K(\mathcal{X})$. For the sample $D = \{(q_i, y_i)\}_{i=1}^n \in (\mathcal{M} \times \mathbb{R})^n$, $x_i = \varphi(q_i) \in \mathbb{R}^p$, the learning gradient algorithm on $\mathcal{M}$ is*

$$\vec{f}_D := \arg\min_{\vec{f} \in \mathcal{H}_K^p} \{\mathcal{E}_D(\vec{f}) + \lambda \|\vec{f}\|_K^2\},$$

*where the weights $w_{ij}$, the RKHS, $\mathcal{H}_K^p$, the RKHS norm $\|\vec{f}\|_K$ and the parameter $\lambda > 0$ are the same as in (2.2).*

Similarly we can deduce the learning gradient algorithm for binary classification on manifolds.

**Definition 2.2.** *Let $\mathcal{M}$ be a Riemannian manifold and $\varphi : \mathcal{M} \to \mathbb{R}^p$ be isometry embedding. Denote $\mathcal{X} = \varphi(\mathcal{M})$ and $\mathcal{H}_K = \mathcal{H}_K(\mathcal{X})$. For the sample $D = \{(q_i, y_i)\}_{i=1}^n \in (\mathcal{M} \times \mathcal{Y})^n$, $x_i = \varphi(q_i) \in \mathbb{R}^p$, the weighted empirical risk for $g : \mathcal{X} \to \mathbb{R}$ and $\vec{f} : \mathcal{X} \to \mathbb{R}^p$*

*is defined as*

$$\mathcal{E}_{\phi,D}(g,\vec{f}) = \sum_{i,j=1}^{n} w_{ij}\phi(y_i(g(x_j) + \vec{f}(x_i) \cdot (x_i - x_j)))$$

*and the algorithm for learning gradient on the manifold is*

$$(g_{\phi,D}, \vec{f}_{\phi,D}) = \arg\min_{(g,\vec{f}) \in \mathcal{H}_K^{p+1}} (\mathcal{E}_{\phi,D}(g,\vec{f}) + \lambda_1\|g\|_K^2 + \lambda_2\|\vec{f}\|_K^2), \tag{2.7}$$

*where $\lambda_1, \lambda_2$ are the regularization parameters.*

Surprisingly these algorithms have forms that are identical to the learning gradient algorithms in Euclidean space. However, the geometric interpretation is different. Note that $\vec{f}$ in Definition 2.1 (or 2.2) models $d\varphi(\nabla_\mathcal{M} f_r)$ (or $d\varphi(\nabla_\mathcal{M} f_c)$), not the gradient itself.

## 3. Convergence rates as a function of the intrinsic dimension

Given the interpretations of the gradient estimates developed in the previous section, it is natural to seek conditions and rates for the convergence of $\vec{f}_D$ to $d\varphi(\nabla_\mathcal{M} f_r)$. Since $I = (d\varphi)^*(d\varphi)$, where $I$ is the identity operator, the convergence to the gradient on the manifold is given by $(d\varphi)^*(f_D) \to \nabla_\mathcal{M} f_r$. The aim of this section is to show that this convergence is true under mild conditions and provide rates.

Throughout this paper we use the following exponential weight function with scale parameter $s^2$,

$$w_{ij} = e^{-\|x_i - x_j\|^2/(2s^2)}.$$

The following $K$-functional will enter our estimates

$$\mathcal{K}(t) = \inf_{\vec{f} \in \mathcal{H}_K^p} (\|(d\varphi)^*\vec{f} - \nabla_\mathcal{M} f_r\|_{L^2_{\rho_\mathcal{M}}}^2 + t\|\vec{f}\|_K^2).$$

The following theorem provides upper bounds for the gradient estimate as a function of the $K$-functional.

**Theorem 3.1.** *Let $\mathcal{M}$ be a compact Riemannian manifold with metric $d_\mathcal{M}$ and let $d\mu$ be the uniform measure on $\mathcal{M}$. Assume the marginal distribution $\rho_\mathcal{M}$ satisfies the regularity conditions:*

(i) *The density $\nu(x) = \frac{d\rho_\mathcal{X}(x)}{d\mu}$ exists and for some $c_1 > 0$ and $0 < \theta \leq 1$*

$$|\nu(x) - \nu(u)| \leq c_1 d_\mathcal{M}^\theta(x,u) \qquad \forall x, u \in \mathcal{M}. \tag{3.1}$$

(ii) *The measure along the boundary is small. There exists $c_2 > 0$ such that*

$$\rho_{\mathcal{M}}(\{x \in \mathcal{M} : d_{\mathcal{M}}(x, \partial \mathcal{M}) \le t\}) \le c_2 t \qquad \forall t > 0. \tag{3.2}$$

*Suppose $f_r \in C^2(\mathcal{M})$. There exists $0 < \varepsilon_0 \le 1$, which depends on $\mathcal{M}$ only, and a constant $C_\rho > 0$ such that, if $s < \varepsilon_0$ and $\lambda = s^{d+2+2\theta}$, then with probability $1 - \delta$ ($\delta \in (0, 1)$) the following bound holds:*

$$\|(d\varphi)^* \vec{f}_D - \nabla_{\mathcal{M}} f_r\|^2_{L^2_{\rho_\mathcal{M}}} \le C_\rho \left(\log \frac{2}{\delta}\right)^2 \left(\frac{1}{n\lambda^2} + s^\theta + \left(\frac{1}{n\lambda^2 s^{2\theta}} + \frac{1}{s^\theta}\right) \mathcal{K}(2s^{2\theta})\right).$$

The rate of convergence of the gradient estimate is an immediate corollary.

**Corollary 3.1.** *Under the assumptions of Theorem 3.1, if $\mathcal{K}(t) = O(t^\beta)$ for some $\frac{1}{2} < \beta \le 1$, then there exist sequences $\lambda = \lambda(n)$ and $s = s(n)$ such that*

$$\|(d\varphi)^* \vec{f}_D - \nabla_{\mathcal{M}} f_r\|^2_{L^2_{\rho_\mathcal{M}}} \to 0 \qquad \text{as } n \to \infty.$$

*If $s = n^{-1/2d+4+5\theta}$ and $\lambda = s^{d+2+2\theta}$, the rate of convergence is*

$$\|(d\varphi)^* \vec{f}_D - \nabla_{\mathcal{M}} f_r\|^2_{L^2_{\rho_\mathcal{M}}} = O(n^{-\theta(2\beta-1)/(2d+4+5\theta)}).$$

This result states that the convergence rate of learning gradient algorithms depends on the intrinsic dimension $d$ of the manifold, not the dimension $p$ of the ambient space. Under the belief that high-dimensional data have low intrinsic dimension $d \ll p$, this explains why the learning gradient algorithms are still efficient for high-dimensional data analysis even when there are limited observations.

If $\mathcal{M}$ is a compact domain in $\mathbb{R}^p$ where $d = p$, $d\varphi = (d\varphi)^* = I$ and $\nabla_{\mathcal{M}}$ is the usual gradient operator, our result reduces to the Euclidean setting that is proven in [16].

The upper bound in Theorem 3.1 is not as tight as possible and may not lead to convergence even when the gradient estimate converges. This is illustrated by the following case. We expect that if $\mathcal{H}_K$ is dense in $C(\mathcal{M})$ then $\mathcal{K}(t) \to 0$ as $t \to 0$ and the gradient estimate should converge in probability to the true gradient. However, Corollary 3.1 states that the convergence holds only when $\mathcal{K}(t)$ decays faster than $O(t^{1/2})$. This is a result of the proof technique we use. More sophisticated but less general proof techniques can give us better estimates and close the above gap; see Remark B.1 in Appendix B for details.

In case of a uniform distribution measure on the manifold, $d\rho_{\mathcal{M}} = d\mu$, we have the following improved upper bound that closes the gap.

**Theorem 3.2.** *Let $\mathcal{M}$ be compact, $d\rho_{\mathcal{M}} = d\mu$ and $f_r \in C^2(\mathcal{M})$. There exists $0 < \varepsilon_0 \le 1$, which depends on $\mathcal{M}$ only, and a constant $C_\mu > 0$ such that, if $s < \varepsilon_0$ and $\lambda = s^{d+3}$, then with probability $1 - \delta$ ($\delta \in (0, 1)$) the following bound holds:*

$$\|(d\varphi)^* \vec{f}_D - \nabla_{\mathcal{M}} f_r\|^2_{L^2_\mu} \le C_\mu \left(\log \frac{2}{\delta}\right)^2 \left(\frac{1}{n\lambda^2} + s + \left(\frac{1}{n\lambda^2 s} + 1\right) \mathcal{K}(2s)\right).$$

An immediate corollary of this theorem is that the rate of convergence of the gradient estimate does not suffer from the same gap as Corollary 3.1.

**Corollary 3.2.** *Under the assumptions of Theorem 3.2, if $\mathscr{K}(t) \to 0$ as $t \to 0$, then*

$$\|(\mathrm{d}\varphi)^* \vec{f}_D - \nabla_{\mathcal{M}} f_r\|^2_{L^2_{\rho_{\mathcal{M}}}} \to 0 \qquad as \ n \to \infty$$

*if $\lambda = s^{d+3}$ and $s = s(n)$ is chosen such that*

$$s \to 0 \quad and \quad \frac{\mathscr{K}(2s^2)}{ns^{2d+7}} \to 0 \qquad as \ n \to \infty.$$

*If in addition $\mathscr{K}(t) = \mathrm{O}(t^\beta)$ for some $0 < \beta \leq 1$, then with the choice $s = n^{-1/(2d+7)}$ and $\lambda = s^{d+3}$, we obtain*

$$\|(\mathrm{d}\varphi)^* \vec{f}_D - \nabla_{\mathcal{M}} f_r\|^2_{L^2_\mu} = \mathrm{O}(n^{-\beta/(2d+7)}).$$

Note that $\mathrm{d}\rho_{\mathcal{M}} = \mathrm{d}\mu$ implies $\nu \equiv 1$ and hence (3.1) holds with $\theta = 1$. In this case the rate in Corollary 3.1 is $\mathrm{O}(n^{-(2\beta-1)/(2d+9)})$. Since $\frac{\beta}{2d+7} > \frac{2\beta-1}{2d+9}$, the rate in Corollary 3.2 is better.

The proofs of the above theorems will be given in Appendix B.

For learning gradients on manifolds for binary classification problems the convergence $(\mathrm{d}\varphi)^* \vec{f}_{\phi,D} \to \nabla_{\mathcal{M}} f_c$ with rate $\mathrm{O}(n^{-1/d})$ can be obtained similarly. We omit the details.

We close this section with some remarks. Note the operator $(\mathrm{d}\varphi)^*$ is the projection onto the tangent space. The convergence results assert that the projection of $\vec{f}_D$ asymptotically approximates the gradient on the manifold. It may be more natural to consider convergence of $\vec{f}_D$ to the gradient on the manifold (after mapping into the ambient space), that is, $\vec{f}_D \to \mathrm{d}\varphi(\nabla_{\mathcal{M}} f_r)$. Unfortunately this is not generally true.

Convergence of learning algorithms that are adaptive to the manifold setting has been considered in the literature [3, 24]. In particular, in [3] local polynomial regression is shown to attain minimax optimal rates in estimating the regression function. Whether it is plausible to extend this to the gradient learning setting is not known. There are a few essential differences between our result and that of Bickel and Li [3]. Unlike our result, their result is pointwise in that convergence and error bounds depend on the point $x \in \mathcal{X}$ and it is not obvious how to obtain $L_2$ convergence from pointwise convergence. In addition, Bickel and Li [3] assume a strong condition on the partial derivatives of the regression function in the ambient space. This may be problematic since these partial derivatives may not be well defined in the ambient space. Since we have different assumptions and a different sense of convergence, the minimax optimality of our results cannot be obtained directly using their arguments. Moreover, in our setting the optimal learning rates will also depend on the choice of the kernel. This is a very interesting open problem.

# 4. Dimension reduction via gradient estimates

The gradient estimates can be used to compute the gradient outer product matrix and provide an estimate of the linear predictive directions or the e.d.r. space.

## 4.1. Estimate gradient outer product matrix and e.d.r. space

Let us start from a semi-parametric model:

$$y = f_r(x) + \epsilon = F(B^T x) + \epsilon,$$

where $B = (b_1, \ldots, b_k) \in \mathbb{R}^{p \times k}$ contains the $k$ e.d.r. directions. When the input space is a domain of $\mathbb{R}^p$ and $f_r$ is smooth, the gradient outer product matrix $G$ with

$$G_{ij} = \left\langle \frac{\partial f_r}{\partial x^i}, \frac{\partial f_r}{\partial x^j} \right\rangle_{L^2}$$

is well defined. It is easy to see that $B$ is given by eigenvectors of $G$ with non-zero eigenvalues.

Using the gradient estimates $\vec{f}_D$ that approximate $\nabla f_r$, we can compute an empirical version of the gradient outer product matrix $\hat{G}$ by

$$\hat{G} = \frac{1}{n} \sum_{\ell=1}^{n} \vec{f}_D(x_\ell) \vec{f}_D(x_\ell)^T.$$

Then the estimates of e.d.r. space can be given by a spectral decomposition of $\hat{G}$.

When the input space is a manifold, since the manifold is not known, we cannot really compute $\nabla_\mathcal{M} f_r$ through $\vec{f}_D$. We can only work directly on $\vec{f}_D$. So we propose to implement the dimension reduction by the same procedure as in the Euclidean setting, that is, we first compute $\hat{G}$ using $\vec{f}_D$ and then estimate the e.d.r. space by the eigenvectors of $\hat{G}$ with top eigenvalues. The only problem is the feasibility of this procedure. The following theorem suggests that the irrelevant dimensions can be filtered out by a spectral decomposition of $\hat{G}$.

**Theorem 4.1.** *Let $\lambda_1 \geq \lambda_2 \geq \cdots \geq \lambda_p$ be the eigenvalues and $u_\ell$, $\ell = 1, \ldots, p$ be the corresponding eigenvectors of $\hat{G}$. Then for any $\ell > k$ we have $B^T u_\ell \to \mathbf{0}_k$ and $u_\ell^T \hat{G} u_\ell \to 0$ in probability.*

## 4.2. An alternative for classification

The idea of dimension reduction by spectral decomposition of the gradient outer production matrix is to estimate the e.d.r. space by finding the directions associated with the directional partial derivative of the largest $L^2$ norm. Intuitively we think the $L^2$ norm may have drawbacks for binary classification problems because a large value of the gradient often is located around the decision boundary, which usually has low density. Thus, an important predictive dimension may correspond to a directional partial derivative with a small $L^2$ norm and hence be filtered out. This motivates us to use $L^\infty$ norm or $\mathcal{H}_K$ norm instead and provide an alternative method.

By using the $\mathcal{H}_K$ norm we consider the gradient covariance matrix (EGCM) $\hat{\Sigma}$ defined by

$$\hat{\Sigma}_{ij} := \langle \vec{f}_{D,i}, \vec{f}_{D,j} \rangle_K \qquad \text{for } 1 \leq i, j \leq p. \tag{4.1}$$

The eigenvectors for the top eigenvalues will be called empirical sensitive linear features (ESF). Estimating the e.d.r. by ESFs is an alternative method for the dimension reduction by gradient estimates and will be referred to as gradient-based linear feature construction (GLFC).

Though using the $L^\infty$ norm or $\mathcal{H}_K$ norm for classification seems to be more intuitive than using the $L^2$ norm, empirically we obtain almost identical results using either method.

### 4.3. Computational considerations

At a first glance one may think it is problematic to compute the spectral decomposition of $\hat{G}$ or $\hat{\Sigma}$ when $p$ is huge. However, due to the special structure of our gradient estimates, they can in fact be realized efficiently in both time, $O(n^2 p + n^3)$, and memory, $O(pn)$. In the following we comment on the computational issues for the EGCM.

In both the regression [16] and the classification [15] settings the gradient estimates satisfy the following representer theorem:

$$\vec{f}_D(x) = \sum_{i=1}^{n} c_{i,D} K(x, x_i).$$

A result of this representer property is that the EGCM has the following positive-semidefinite quadratic form

$$\hat{\Sigma} = c_D K c_D^T,$$

where $c_D = (c_{1,D}, \ldots, c_{n,D}) \in \mathbb{R}^{p \times n}$ and $K$ is the $n \times n$ kernel matrix with $K_{ij} = K(x_i, x_j)$. Since $K$ is positive definite, there exists a lower-diagonal matrix $L$ so that $LL^T = K$. Then

$$\hat{\Sigma} = \tilde{c}_D \tilde{c}_D^T,$$

where $\tilde{c}_D = c_D L$ is $p \times n$ matrix. An immediate result of this formula is that the EGCM is a rank $n$ matrix and has at most $n$ non-zero eigenvalues and therefore at most the top $n$ empirical features will be selected as relevant ones. Efficient solvers for the first $n$ eigenvectors of $\hat{\Sigma}$ using QR decomposition for $\tilde{c}_D$ [10] are standard and require $O(n^2 p)$ time and $O(pn)$ memory. This observation conforms to the intuition that with $n$ samples, it is impossible to select more than $n$ independent features or predict a function that depends on more than $n$ independent variables.

## 4.4. Relations to OPG method

MAVE and the OPG method proposed in [23] share similar ideas by using the gradient outer product matrix implicitly or explicitly. Of them, OPG is more related to learning gradients. We discuss differences between the methods.

At each point $x_j$ the OPG method estimates the function value $a_j \approx f_r(x_j)$ and gradient vector $b_j \approx \nabla f_r(x_j)$ by

$$(a_j, b_j) = \arg\min_{a \in \mathbb{R}, b \in \mathbb{R}^p} \sum_{i=1}^n w_{ij}[y_i - a - b^\mathrm{T}(x_i - x_j)]^2. \quad (4.2)$$

Then the gradient outer product matrix is approximated by

$$\hat{G} = \frac{1}{n} \sum_{j=1}^n b_j b_j^\mathrm{T}.$$

Notice that if we set the kernel as $K(x, u) = \delta_{x,u}$ and $\lambda = 0$ the learning gradient algorithm reduces to (4.2). In this sense the learning gradient extends OPG from estimating the gradient vector only at the sampling points to estimating the gradients by a vector valued function. Hence the estimates extend to out-of-sample points. This offers the potential to apply the method more generally; for example, numerical derivatives in a low-dimensional space or function adaptive diffusion maps (see [17] for details).

The solution to the minimization problem in (4.2) is not unique when $p > n$. This can result in overfitting and instability of the estimate of the gradient outer product matrix. In this sense OPG cannot be directly used for the "large $p$, small $n$" problem. The regularization term in the learning gradient algorithms helps to reduce overfitting and allows for feasible estimates in the "large $p$, small $n$" setting. However, we should remark that this is a theoretical and conceptual argument. In practice OPG can be used together with preprocessing the data using principal components analysis (PCA) and results in performance comparable to learning gradients.

## 5. Results on simulated and real data

In this section we illustrate the utility and properties of our method. We will focus on the "large $p$, small $n$" setting and binary classification problems.

In the following simulations we always set the weight function as $\exp(-\frac{\|x-u\|^2}{s^2})$ with $s$ as the median of pairwise distance of the sampling points. The Mercer kernel $K$ is $K(x, u) = \exp(-\frac{\|x-u\|^2}{\sigma^2})$ with $\sigma$ equal to 0.2 times the median of pairwise distance of the sample points. They may not be optimal but work well when $p \gg n$ in our experience.

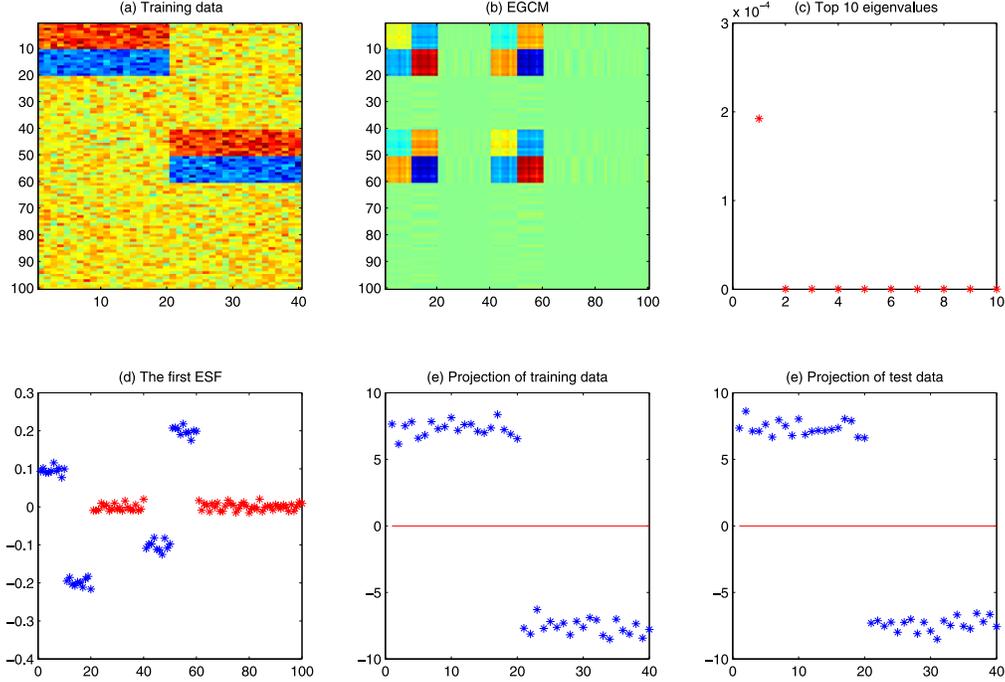

**Figure 1.** Linear classification simulation with $\sigma = 0.5$.

## 5.1. A linear classification simulation

Data are drawn from two classes in an $n = 100$ dimensional space. Samples from class $-1$ were drawn from

$$x^{j=1:10} \sim \mathrm{N}(1.5, \sigma), \qquad x^{j=11:20} \sim \mathrm{N}(-3, \sigma), \qquad x^{j=21:100} \sim \mathrm{N}(0, \sigma).$$

Samples from class $+1$ were drawn from

$$x^{j=41:50} \sim \mathrm{N}(-1.5, \sigma), \qquad x^{j=51:60} \sim \mathrm{N}(3, \sigma), \qquad x^{j=1:40, 61:100} \sim \mathrm{N}(0, \sigma).$$

Note that $\sigma$ measures the noise level and difficulty of extracting the correct dimensions. We drew 20 observations for each class from the above distribution as training data and another 40 samples are independently drawn as test data. By changing the noise level $\sigma$ from 0.2 to 3, we found our method stably finds the correct predictive directions when $\sigma < 2$. From a prediction point of view the result is still acceptable when $\sigma > 2$ though the estimates of the predictive dimension contain somewhat larger errors. In Figures 1 and 2 we show the results for $\sigma = 0.5$ and $\sigma = 2.5$, respectively.

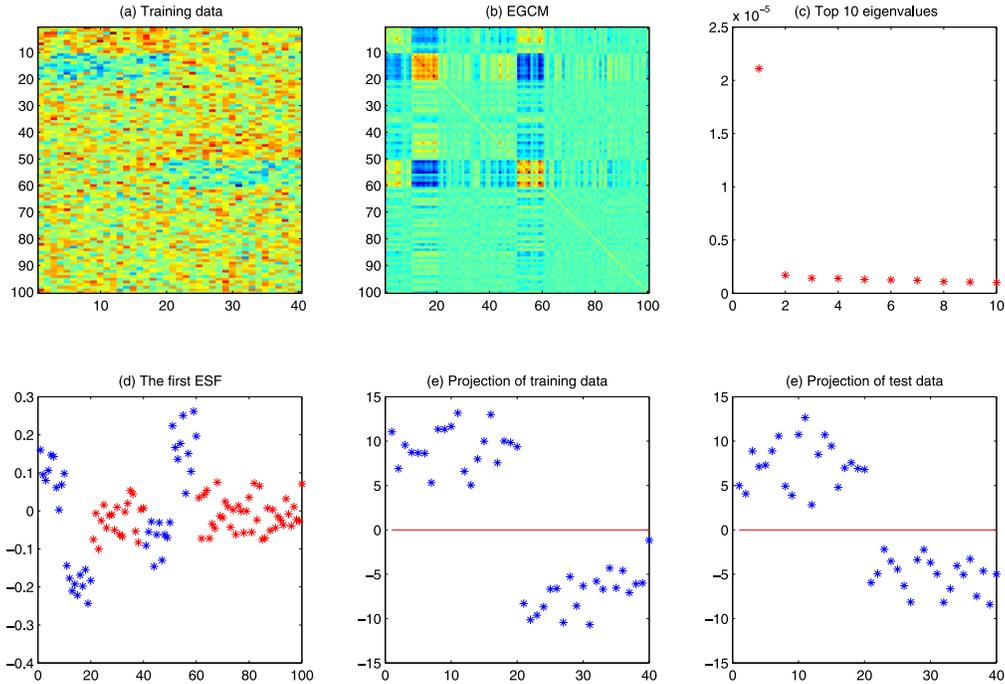

**Figure 2.** Linear classification simulation with $\sigma = 2.5$.

## 5.2. A nonlinear classification simulation

Data are drawn from two classes in a $p$-dimensional space. Only the first $d$-dimensions are relevant in the classification problem. For samples from class $+1$ the first $d$-dimensions correspond to points drawn uniformly from the surface of a $d$-dimensional hypersphere with radius $r$, for class $-1$ the first 10 dimensions correspond to points drawn uniformly from the surface of a $d$-dimensional hypersphere with radius $2.5r$. The remaining $p$-$d$-dimensions are noise

$$x^j \sim \mathrm{N}(0, \sigma) \qquad \text{for } j = d+1, \ldots, p.$$

Note that the data can be separated by a hypersphere in the first $d$-dimensions. Therefore projecting the data onto the first $d$-ESFs for this problem should reflect the underlying geometry.

In this simulation we set $p = 200$, $d = 2$, $r = 3$ and $\sigma$ varying from 0.1 to 1. The first two ESFs are shown to capture the correct underlying structure when $\sigma \leq 0.7$. In Figure 3 we give the result with $\sigma = 0.2$. We also studied the influence of $p$ and found when $p \leq 50$ the noise level can be as large as 1.0.

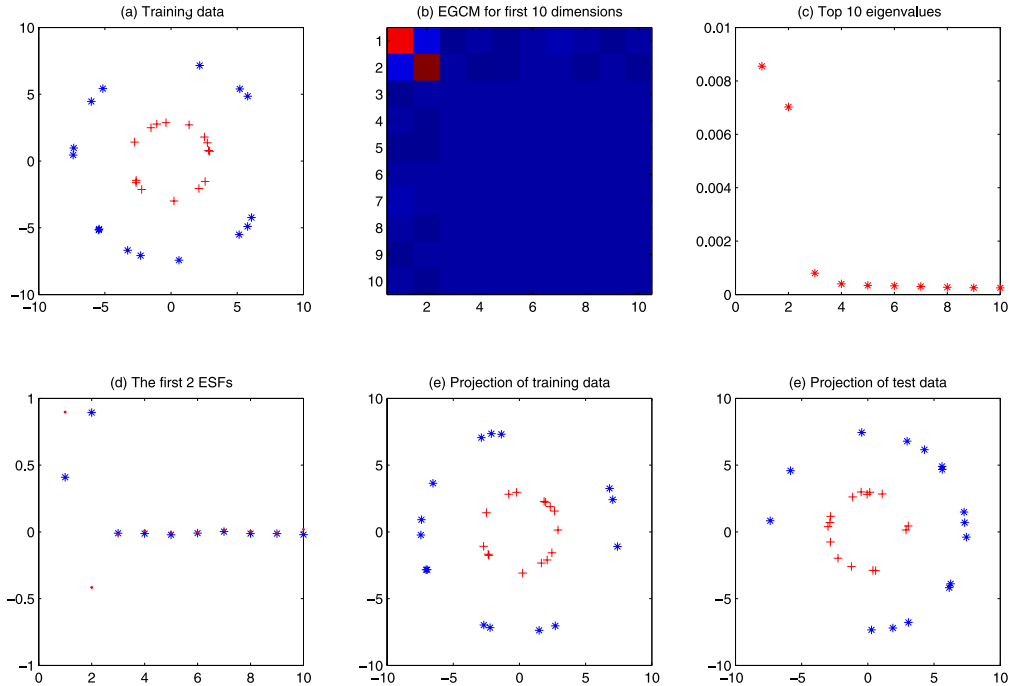

**Figure 3.** Nonlinear classification simulation with $\sigma = 0.2$. In (d) the first ESF is in blue and the second in red.

## 5.3. Digit classification

A standard data set used in the machine learning community to benchmark classification algorithms is the MNIST data set (Y. LeCun, http://yann.lecun.com/exdb/mnist/). The data set contains 60 000 images of handwritten digits $\{0, 1, 2, \ldots, 9\}$, where each image consists of $p = 28 \times 28 = 784$ grayscale pixel intensities. This data set is commonly believed to have strong nonlinear manifold structure. In this section we report results on one of the most difficult pairwise comparisons: discriminating a handwritten "3" from an "8".

In the simulation we randomly choose 30 images from each class and the remaining are used as the test set. We compare the following dimension reduction methods GLFC, SIR and OPG. In Table 1 we report the classification error rates by the k-NN classifier with $k = 5$ using the respective method for dimension reduction. The SIR results reported are for a regularized version of SIR (RSIR) since SIR is not stable for very high-dimensional data. As mentioned before OPG cannot be directly applied so we first run PCA. We compare the results for using all PCs (PC-OPG) and 30 PCs (PC30-OPG), respectively. The last column is the error rate by k-NN without dimension reduction.

### 5.4. Gene expression data

One problem domain where high dimensions are ubiquitous is the analysis and classification of gene expression data. We consider two classification problems based on gene expression data. One is the study using expression data to discriminate acute myeloid leukemia (AML) from acute lymphoblastic leukemia (ALL) [11] and another is the classification of prostate cancer [19]. In the leukemia data set there are 48 samples of AML and 25 samples of ALL. The number of genes is $p = 7129$. The data set was split into a training set of 38 samples and a test set of 35 samples as specified in [11]. In the prostate cancer data, the dimension is $p = 12\,600$. The training data contains 102 samples, 52 tumor samples and 50 non-tumor samples. The independent test data contains 34 samples from a different experiment. We applied GLFC, SIR and OPG to these two data sets and compared the accuracy using a linear support vector machine classifier. The leave-one-out (LOO) error over the training data and the test error are reported in Tables 2 and 3 for the leukemia data and prostate cancer data, respectively. For leukemia classification the two classes are well separated and all methods perform quite similarly. For prostate cancer classification, the accuracy in [19] is about 90%. All methods achieve similar accuracy on the training data. GLFC and OPG methods have better prediction power on test data. From our experiments (both for gene expression data and digits data) we see that the PC-OPG method performs quite similarly to GLFC when the number of top PCs is correctly set. But it seems quite sensitive to this choice.

## 6. Discussion

In this paper we first extended the gradient estimation and feature selection framework outlined in [15, 16] from the ambient space setting to the manifold setting. Convergence is shown to depend on the intrinsic dimension of the manifold but not the dimension of the ambient Euclidean space. This helps to explain the feasibility "large

**Table 1.** Classification error rate for 3 vs. 8

| GLFC | PC-OPG | PC30-OPG | SIR | RSIR | kNN ($k = 5$) |
|---|---|---|---|---|---|
| 0.0853 | 0.1110 | 0.0912 | 0.1877 | 0.0956 | 0.1024 |

**Table 2.** Classification error for leukemia data

|  | GLFC | PC-OPG | SIR | SVM |
|---|---|---|---|---|
| LOO error | 1 | 1 | 1 | 1 |
| Test error | 0 | 0 | 1 | 1 |
| Dimension | $d = 2$ | $d = 4$ | $d = 1$ | $d = 7129$ |

**Table 3.** Classification error for prostate cancer data

|  | GLFC | PC-OPG | PC50-OPG | SIR | SVM |
|---|---|---|---|---|---|
| LOO error | 9 | 15 | 10 | 9 | 9 |
| Test error | 2 | 9 | 6 | 9 | 9 |
| Dimension | $d=5$ | $d=2$ | $d=1$ | $d=1$ | $d=12\,600$ |

$p$, small $n$" problems. We outlined properties of this method and illustrated its utility for real and simulated data. Matlab code for learning gradients can be obtained at http://www.stat.duke.edu/~sayan/soft.html.

We close by stating open problems and discussion points:

1. Large $p$, not so small $n$: The computational approaches used in this paper require that $n$ is small. The theory we provide does not place any constraint on $n$. The extension of the computational methods to larger $n$ involves the ability to expand the gradient estimates in terms of an efficient bases expansion. For the approach proposed in this paper the number bases are at most $n^2$, which is efficient for small $n$.
2. Fully Bayesian model: The Tikhonov regularization framework coupled with the use of an RKHS allows us to implement a fully Bayesian version of the procedure in the context of Bayesian radial basis (RB) models [14]. The Bayesian RB framework can be extended to develop a proper probability model for the gradient learning problem. The optimization procedures in Definitions 2.1 and 2.2 would be replaced by Markov chain Monte Carlo methods and the full posterior rather than the maximum a posteriori estimate would be computed. A very useful result of this is that in addition to the point estimates for the gradient we would also be able to compute confidence intervals.

## Appendix A: Geometric background

In this section, we introduce some properties on manifolds that we will need in our proofs. Let $\mathcal{M}$ be a Riemannian manifold and $\varphi:\mathcal{M}\to\mathbb{R}^p$ be an isometric embedding, that is, for every $q\in\mathcal{M}$, $\mathrm{d}\varphi_q:T_q\mathcal{M}\to T_{\varphi(q)}\mathbb{R}^p\cong\mathbb{R}^p$ is isometric:

$$\mathrm{d}\varphi_q(v)\cdot\mathrm{d}\varphi_q(v_1)=\langle v,v_1\rangle_q \qquad \forall v,v_1\in T_q\mathcal{M}.$$

This directly leads to the following conclusion.

**Lemma A.1.** *For every $q\in\mathcal{M}$, the following hold:*

(i) $(\mathrm{d}\varphi_q)^*\circ(\mathrm{d}\varphi_q)=I_{T_q\mathcal{M}}$, *the identity operator on $T_q\mathcal{M}$;*
(ii) $(\mathrm{d}\varphi_q)\circ(\mathrm{d}\varphi_q)^*$ *is the projection operator of $T_{\varphi(q)}\mathbb{R}^p$ to its subspace $\mathrm{d}\varphi_q(T_q\mathcal{M})\subset T_{\varphi(q)}\mathbb{R}^p$;*

(iii) $\|\mathrm{d}\varphi_q\| = \|(\mathrm{d}\varphi_q)^*\| = 1$.

**Lemma A.2.** *Let $\mathcal{M}$ be compact. There exists $\varepsilon_0 > 0$ uniform in $q \in \mathcal{M}$ such that $\exp_q$ is well defined on $B_q(\epsilon_0)$ and is a diffeomorphism onto its image. Moreover, given an orthonormal basis $\{e_i^q\}_{i=1}^d$ of $T_q\mathcal{M}$, under the q-normal coordinates defined by $\exp_q^{-1}$, if $|v| \leq \varepsilon_0$ the following hold:*

(i) $\frac{1}{2} \leq \sqrt{\det(G)}(v) \leq \frac{3}{2}$ *where $G$ is the volume element.*
(ii) $\frac{1}{2}|v|^2 \leq \|\varphi(\exp_q(v)) - \varphi(q)\| \leq |v|^2$.
(iii) $\varphi(\exp_q(v)) - \varphi(q) = \mathrm{d}\varphi_q(v) + \mathrm{O}(|v|^2)$.

**Proof.** This follows directly from the compactness of $\mathcal{M}$ and Proposition 2.2 in [9]. See [24] for a self-contained proof of a very similar result. □

**Lemma A.3.** *Let $\mathcal{M}$ be compact and $\varepsilon_0$ be given as in Lemma A.2. If $f \in C^2(\mathcal{M})$, then there exists a constant $C > 0$ such that for all $q \in \mathcal{M}$ and $v \in T_q\mathcal{M}$, $|v| \leq \varepsilon_0$,*

$$|f(\exp_q(v)) - f(q) - \langle \nabla_\mathcal{M} f(q), v \rangle| \leq C|v|^2.$$

**Proof.** Since $f \in C^2(\mathcal{M})$, $f \circ \exp_q(v)$ is $C^2$ on $B_q(\varepsilon_0)$. By the discussion in Section 2.2,

$$\begin{aligned}
&|f(\exp_q(v)) - f(q) - \langle \nabla_\mathcal{M} f(q), v \rangle| \\
&= |(f \circ \exp_q)(v) - (f \circ \exp_q)(0) - \nabla(f \circ \exp_q)(0) \cdot v| \\
&\leq C_q |v|^2
\end{aligned}$$

with $C_q = \sup_{v \in B_q(\varepsilon_0)} |\nabla^2 (f \circ \exp_q)(v)|$. Since $\nabla^2(f \circ \exp_q)(v)$ is continuous in $q$ and $\mathcal{M}$ is compact, $C = \sup_{q \in \mathcal{M}} C_q$ exists and our conclusion follows. □

We remark that if $\mathcal{M} \subset \mathbb{R}^p$ is a submanifold with intrinsic metric, $I_{\mathbb{R}^p}$ is an isometric embedding. If in addition $\mathcal{M}$ is a closed domain in $\mathbb{R}^p$, then $\exp_x(v) = x + v$ for $x \in \mathcal{M}$.

## Appendix B: Proofs of convergence results

In the manifold setting we usually do not know the manifold or the embedding. It will be convenient to regard $\mathcal{M}$ as a submanifold of $\mathbb{R}^p$ and $d_\mathcal{M}$ as the intrinsic metric on the manifold. This implies that $\varphi = I_{\mathbb{R}^p}$ and we can identify $\mathcal{M}$ as $\mathcal{X} = \varphi(\mathcal{M})$. Note that the marginal distribution $\rho_\mathcal{M}$ on $\mathcal{M}$ induces a distribution $\rho_\mathcal{X}$ on $\mathcal{X} = \varphi(\mathcal{M})$. The above notation implies $\rho_\mathcal{M} = \rho_\mathcal{X}$. We denote $\mathscr{D}_x = \mathrm{d}(I_{\mathbb{R}^p})_x$ and $\mathscr{D}$ is the operator on vector fields such that $\mathscr{D}h(x) = \mathscr{D}_x(h(x))$ for $h \in T\mathcal{M}$. Correspondingly, $\mathscr{D}_x^*$ is the dual of $\mathscr{D}_x$ and $\mathscr{D}^*$ maps a $p$-dimensional vector valued function $\vec{f}$ to a vector field $\mathscr{D}^*\vec{f}(x) = \mathscr{D}_x^*(\vec{f}(x))$. We will adopt this notation to simplify our expression and give the proofs of the results in Section 3.

Recall the following definition given in [16].

**Definition B.1.** *Denote $Z = X \times Y = \mathcal{M} \times Y$. For $s > 0$ and $\vec{f} : \mathcal{M} \to \mathbb{R}^p$, define the expected error*

$$\mathcal{E}(\vec{f}) = \int_Z \int_Z e^{-\|x-\xi\|^2/(2s^2)} (y - \eta + \vec{f}(x) \cdot (\xi - x))^2 \, d\rho(x,y) \, d\rho(\xi,\eta).$$

*If we denote $\sigma_s^2 = \int_Z \int_Z e^{-\|x-\xi\|^2/(2s^2)} (y - f_r(x))^2 \, d\rho(x,y) \, d\rho(\xi,\eta)$, then*

$$\mathcal{E}(\vec{f}) = 2\sigma_s^2 + \int_{\mathcal{M}} \int_{\mathcal{M}} e^{-\|x-\xi\|^2/(2s^2)} (f_r(x) - f_r(\xi) + \vec{f}(x) \cdot (\xi - x))^2 \, d\rho_{\mathcal{M}}(x) \, d\rho_{\mathcal{M}}(\xi).$$

Define

$$\vec{f}_\lambda = \arg \min_{\vec{f} \in \mathcal{H}_K^p} \{\mathcal{E}(\vec{f}) + \lambda \|\vec{f}\|_K^2\}.$$

It can be regarded as the infinite sample limit of $\vec{f}_{D,\lambda}$. The following decomposition holds

$$\|\mathscr{D}^* \vec{f}_D - \nabla_{\mathcal{M}} f_r\|_{L^2_{\rho_{\mathcal{M}}}} \leq \kappa \|\vec{f}_D - \vec{f}_\lambda\|_K + \|\mathscr{D}^* \vec{f}_\lambda - \nabla_M f_r\|_{L^2_{\rho_{\mathcal{M}}}}.$$

This bounds $\|\mathscr{D}^* \vec{f}_D - \nabla_{\mathcal{M}} f_r\|_{L^2_{\rho_{\mathcal{M}}}}$ by two terms. The first one is called *sample error* and the second is the *approximation error*.

For the sample error, we have the following estimate.

**Proposition B.1.** *Assume $|y| \leq M$ almost surely. There are two constants $C_1, C_2 > 0$ such that for any $\delta > 0$, with confidence $1 - \delta$,*

$$\|\vec{f}_D - \vec{f}_\lambda\|_K \leq \frac{C_1 \log(2/\delta)}{\sqrt{m}\lambda} + \left( \frac{C_2 \log(2/\delta)}{\sqrt{m}\lambda} + \frac{1}{m} \right) \|\vec{f}_\lambda\|_K.$$

This estimate has in fact been given in [16]. To see this, we notice two facts. First, our algorithm in Definition 2.1 is different from that in [16] only by a scalar. Second, the proof in [16] does not depend on the geometric structure of the input space. So a scalar argument leads to the above bound for the sample error directly. Of course, one may obtain better estimates by incorporating the specific geometric property of the manifold.

Next we turn to estimate the approximation error. In the following, we always assume $\mathcal{M}$ compact and set $\varepsilon_0$ to be the same as in Lemma A.2. Without loss of generality, we also assume $\varepsilon_0 \leq 1$.

**Proposition B.2.** *Assume (3.1) and (3.2). If $f_r \in C^2(\mathcal{M})$, there is a constant $C_3 > 0$ such that for all $\lambda > 0$ and $s < \varepsilon_0$,*

$$\|\mathscr{D}^* \vec{f}_\lambda - \nabla_{\mathcal{M}} f_r\|_{L^2_{\rho_{\mathcal{M}}}}^2 \leq C_3 \left( s^\theta + \left( \frac{s^{d+2+\theta}}{\lambda} + \frac{1}{s^\theta} \right) \left( s^2 + \mathscr{K} \left( \frac{\lambda}{s^{d+2}} + s^2 \right) \right) \right).$$

*If $d\rho_{\mathcal{M}} = d\mu$, the estimate can be improved.*

**Proposition B.3.** *Let $f_r \in C^2(\mathcal{M})$. If $d\rho_\mathcal{M} = d\mu$, then there exists a constant $C_{3,\mu} > 0$ such that for all $\lambda > 0$ and $s < \varepsilon_0$,*

$$\|\mathscr{D}^*\vec{f}_\lambda - \nabla_\mathcal{M} f_r\|_{L_\mu^2}^2 \leq C_{3,\mu}\left(s + \left(\frac{s^{d+3}}{\lambda} + 1\right)\left(s^2 + \mathscr{K}\left(\frac{\lambda}{s^{d+2}} + s^2\right)\right)\right).$$

The proof of these bounds for the approximation error will be given in two steps. In the first step we bound the $L^2$-difference by the expected error and in the second step the functional $\mathscr{K}$ is used to control the expected error.

**Lemma B.1.** *Assume Condition 3.2 and $f_r \in C^2(\mathcal{M})$. There exists a constant $c_3 > 0$ so that for all $s < \varepsilon_0$ and $\vec{f} \in \mathcal{H}_K^n$,*

$$\|\mathscr{D}^*\vec{f} - \nabla_\mathcal{M} f\|_{L_{\rho_\mathcal{M}}^2}^2 \leq c_3\left((1 + \|\vec{f}\|_K)^2 s^\theta + \frac{1}{s^{d+2+\theta}}(\mathcal{E}(\vec{f}) - 2\sigma^2)\right). \tag{B.1}$$

*If, in addition, $d\rho_\mathcal{M} = d\mu$, then the estimate can be improved to*

$$\|\mathscr{D}^*\vec{f} - \nabla_\mathcal{M} f\|_{L_{\rho_\mathcal{M}}^2}^2 \leq c_{3,\mu}\left((1 + \|\vec{f}\|_K)^2 s + \frac{1}{s^{d+2}}(\mathcal{E}(\vec{f}) - 2\sigma^2)\right) \tag{B.2}$$

*for some $c_{3,\mu} > 0$.*

**Proof.** Denote $X_s = \{x \in \mathcal{M} : d_\mathcal{M}(x, \partial M) \geq s$ and $\nu(x) \geq (1 + c_2)s^\theta\}$. For $x \in \mathcal{M}$, let $\mathcal{B}_{x,s} = \{\xi \in \mathcal{M} : d_\mathcal{M}(\xi, x) \leq s\}$. We prove the conclusion in three steps.

*Step* 1. Define the local error function

$$\text{er}_s(x) = \int_{\mathcal{B}_{x,s}} e^{-\|x-\xi\|^2/(2s^2)}(f_r(x) - f_r(\xi) + \vec{f}(x) \cdot (\xi - x))^2 \, d\mu(\xi).$$

We claim that there exists a constant $c' > 0$ such that for each $x \in X_s$,

$$|\mathscr{D}^*\vec{f}(x) - \nabla_\mathcal{M} f_r| \leq c'\left((1 + \|\vec{f}\|_K)^2 s^2 + \frac{1}{s^{d+2}} \text{er}_s(x)\right). \tag{B.3}$$

Since $s < \varepsilon_0$, $\exp_x$ is a homeomorphism from $B_x(s)$ onto $\mathcal{B}_{x,s}$. For every $\xi \in \mathcal{B}_{x,s}$ there exists $v \in B_x(s)$ so that $\xi = \exp_x(v)$. Write every $v \in T_x\mathcal{M}$ in normal coordinates. Then $\text{er}_s(x)$ equals

$$\int_{B_x(s)} e^{-\|x-\exp_x(v)\|^2/(2s^2)}(f_r(x) - f_r(\exp_x(v)) + \vec{f}(x) \cdot (\exp_x(v) - x))^2 \sqrt{\det(G)}(v) \, dv.$$

Denote

$$t_1(v) = f_r(x) - f_r(\exp_x(v)) + \vec{f}(x) \cdot (\exp_x(v) - x) - \langle \mathscr{D}^*\vec{f}(x) - \nabla_\mathcal{M} f_r(x), v\rangle,$$
$$t_2(v) = \langle \mathscr{D}^*\vec{f}(x) - \nabla_\mathcal{M} f_r(x), v\rangle.$$

By the assumption $f_r \in C^2(\mathcal{M})$ and Lemma A.3,

$$|f_r(x) - f_r(\exp_x(v)) + \langle \nabla_{\mathcal{M}} f_r(x), v \rangle| \leq \tilde{c}_1 |v|^2$$

for some $\tilde{c}_1 > 0$. By Lemma A.2(iii), there exists $\tilde{c}_2 > 0$ depending on $\mathcal{M}$ only so that

$$|\vec{f}(x) \cdot (\exp_x(v) - x) - \vec{f}(x) \cdot \mathscr{D}_x(v)| \leq \tilde{c}_2 \|\vec{f}(x)\| |v|^2 \leq \tilde{c}_2 \kappa \|\vec{f}\|_K |v|^2.$$

Notice that $\vec{f}(x) \cdot \mathscr{D}_x(v) = \langle \mathscr{D}_x^*(\vec{f}(x)), v \rangle = \langle \mathscr{D}^* \vec{f}(x), v \rangle$. So we have

$$|t_1(v)| \leq (\tilde{c}_1 + \tilde{c}_2 \kappa \|\vec{f}\|_K) |v|^2.$$

Using the facts $\frac{1}{2}|v|^2 \leq \|x - \exp_x(v)\|^2 \leq |v|^2$ and $\frac{1}{2} \leq \sqrt{\det(G)}(v) < \frac{3}{2}$, we obtain

$$\int_{B_x(s)} e^{-\|x - \exp_x(v)\|^2/(2s^2)} |t_1(v)|^2 \sqrt{\det(G)}(v) \, dv$$

$$\leq \frac{3}{2} (\tilde{c}_1 + \tilde{c}_2 \|\vec{f}\|_K)^2 \int_{|v| \leq s} e^{-|v|^2/(4s^2)} |v|^4 \, dv$$

$$\leq \frac{3}{2} (\tilde{c}_1 + \tilde{c}_2 \|\vec{f}\|_K)^2 \tilde{c}_3 s^{d+4},$$

where $\tilde{c}_3 = \int_{|v| \leq 1} e^{-|v|^2/4} |v|^4 \, dv$.

Denote

$$Q_s(x) = \int_{B_x(s)} e^{-\|x - \exp_x(v)\|^2/(2s^2)} |t_2(v)|^2 \sqrt{\det(G)}(v) \, dv.$$

By the Schwarz inequality

$$\int_{B_x(s)} e^{-\|x - \exp_x(v)\|^2/(2s^2)} |t_1(v)| |t_2(v)| \sqrt{\det(G)}(v) \, dv$$

$$\leq \left( \int_{B_x(s)} e^{-\|x - \exp_x(v)\|^2/(2s^2)} |t_1(v)|^2 \sqrt{\det(G)}(v) \, dv \right)^{1/2} \sqrt{Q_s(x)}$$

$$\leq (\tilde{c}_1 + \tilde{c}_2 \|\vec{f}\|_K) \sqrt{\frac{3}{2} \tilde{c}_3 s^{d+4}} \sqrt{Q_s(x)}.$$

Then we get

$$\mathrm{er}_s(x) \geq \int_{B_x(s)} e^{-\|x - \exp_x(v)\|^2/(2s^2)} (|t_2(v)|^2 - 2|t_1(v)||t_2(v)| - 2|t_1(v)|^2) \sqrt{\det(G)}(v) \, dv$$

$$\geq Q_s(x) - (\tilde{c}_1 + \tilde{c}_2 \|\vec{f}\|_K) \sqrt{\frac{3}{2} \tilde{c}_3 s^{d+4}} \sqrt{Q_s(x)} - \frac{3}{2} (\tilde{c}_1 + \tilde{c}_2 \|\vec{f}\|_K)^2 \tilde{c}_3 s^{d+4}.$$

This implies

$$Q_s(x) \leq 3(\tilde{c}_1 + \tilde{c}_2 \|\vec{f}\|_K)^2 \tilde{c}_3 s^{d+4} + 2\operatorname{er}_s(x).$$

By the facts $\|x - u\|^2 \leq |v|^2$, $\sqrt{\det(G)}(v) \geq \frac{1}{2}$ and $\int_{B_x(s)} e^{-|v|^2/(2s^2)} v_i v_j = 0$, if $i \neq j$, we obtain

$$Q_s(x) \geq \frac{1}{2} \sum_{i,j=1}^{d} (\mathscr{D}^* \vec{f}(x) - \nabla_\mathcal{M} f_r(x))_i (\mathscr{D}^* \vec{f}(x) - \nabla_\mathcal{M} f_r(x))_j \int_{B_x(s)} e^{-|v|^2/(2s^2)} v_i v_j \, dv$$

$$= \tilde{c}_4 s^{d+2} |\mathscr{D}^* \vec{f}(x) - \nabla_\mathcal{M} f_r(x)|^2,$$

where $\tilde{c}_4 = \frac{1}{2d} \int_{|v| \leq 1} e^{-|v|^2/2} |v|^2 \, dv$. Therefore, our claim (B.3) holds with

$$c' = \frac{1}{\tilde{c}_4}(3(\tilde{c}_1 + \tilde{c}_2)^2 \tilde{c}_3 + 2).$$

*Step* 2. By (B.3) we have

$$\int_{X_s} |\mathscr{D}^* \vec{f}(x) - \nabla_\mathcal{M} f_r(x)|^2 \, d\rho_\mathcal{M}(x) \tag{B.4}$$
$$\leq c' \left( (1 + \|\vec{f}\|_K)^2 s^2 + \frac{1}{s^{d+2}} \int_{X_s} \operatorname{er}_s(x) \, d\rho_\mathcal{M}(x) \right).$$

By the assumption $d\rho_\mathcal{M}(\xi) = \nu(\xi) \, d\mu$ and (3.2), we have $\nu(\xi) \geq s^\theta$ if $x \in X_s$ and $\xi \in \mathcal{B}_{x,s}$. Therefore,

$$\int_{\mathcal{B}_{x,s}} e^{-\|x - \xi\|^2/2s^2} (f_r(x) - f_r(\xi) + \vec{f}(x) \cdot (\xi - x))^2 \, d\rho_\mathcal{M}(\xi) \geq s^\theta \operatorname{er}_s(x).$$

Integrating both sides over $x$ on $X_s$ and using the fact $\mathcal{B}_{x,s} \subset \mathcal{M}$ when $x \in X_s$, we obtain

$$\int_{X_s} \operatorname{er}_s(x) \, d\rho_\mathcal{M}(x) \leq \frac{1}{s^\theta}(\mathcal{E}(\vec{f}) - 2\sigma_s^2).$$

Plugging into (B.4) gives

$$\int_{X_s} |\mathscr{D}^* \vec{f}(x) - \nabla_\mathcal{M} f_r(x)|^2 \, d\rho_\mathcal{M}(x) \leq c' \left( (1 + \|\vec{f}\|_K)^2 s^2 + \frac{1}{s^{d+2+\theta}}(\mathcal{E}(\vec{f}) - 2\sigma_s^2) \right). \tag{B.5}$$

If $d\rho_\mathcal{M} = d\mu$, we immediately obtain

$$\int_{X_s} \operatorname{er}_s(x) \, d\rho_\mathcal{M}(x) \leq \mathcal{E}(\vec{f}) - 2\sigma_s^2$$

and hence

$$\int_{X_s} |\mathscr{D}^* \vec{f}(x) - \nabla_{\mathcal{M}} f_r(x)|^2 \, \mathrm{d}\rho_{\mathcal{M}}(x) \leq c' \bigg( (1 + \|\vec{f}\|_K)^2 s^2 + \frac{1}{s^{d+2}} (\mathcal{E}(\vec{f}) - 2\sigma_s^2) \bigg). \quad \text{(B.6)}$$

*Step* 3. Condition (3.1) implies $\nu$ is continuous. Since $\mathcal{M}$ is compact, $\sup_{x \in \mathcal{M}} \nu(x) = \tilde{c}_5$ exists. So

$$\rho_{\mathcal{M}}(\mathcal{M} \setminus X_s) \leq c_2 s + \tilde{c}_5 (1 + c_2) s^\theta \mu(\mathcal{M}).$$

Together with the fact $|\mathscr{D}^* \vec{f}(x) - \nabla_{\mathcal{M}} f_r(x)| \leq \kappa \|\vec{f}\|_K + \|\nabla_{\mathcal{M}} f_r\|_\infty$, we have

$$\int_{\mathcal{M} \setminus X_s} |\mathscr{D}^* \vec{f}(x) - \nabla_{\mathcal{M}} f_r(x)|^2 \, \mathrm{d}\rho_{\mathcal{M}}(x) \leq \tilde{c}_6 (1 + \|\vec{f}\|_K)^2 s^\theta \quad \text{(B.7)}$$

with $\tilde{c}_6 = (\kappa + \|\nabla_{\mathcal{M}} f_r\|_\infty)(c_2 + \tilde{c}_5 (1 + c_2) \mu(\mathcal{M}))$.

Combining (B.5) and (B.7) leads to conclusion (B.1).

If $\mathrm{d}\rho_{\mathcal{M}} = \mathrm{d}\mu$, $\nu(x) = 1$ for all $x \in \mathcal{M}$, (3.1) holds with $\theta = 1$ and $\tilde{c}_5 = 1$. So (B.7) holds with $\theta = 1$. This together with (B.6) proves (B.2).

This finishes the proof. □

Define the functional

$$\mathcal{A}(s, \lambda) = \inf_{\vec{f} \in \mathcal{H}_K^n} (\mathcal{E}(\vec{f}) - 2\sigma_s^2 + \lambda \|\vec{f}\|_K^2).$$

Applying Lemma B.1 to $\vec{f}_\lambda$, we immediately obtain the following corollary.

**Corollary B.1.** *Under Assumption (B.1), we have*

$$\|\mathscr{D}^* \vec{f}_\lambda - \nabla_{\mathcal{M}} f_r\|_{L^2_{\rho_{\mathcal{M}}}}^2 \leq c_3 \bigg( 2s^\theta + \bigg( \frac{2s^\theta}{\lambda} + \frac{1}{s^{d+2+\theta}} \bigg) \mathcal{A}(s, \lambda) \bigg).$$

*If* $\mathrm{d}\rho_{\mathcal{M}} = \mathrm{d}\mu$, *then*

$$\|\mathscr{D}^* \vec{f}_\lambda - \nabla_{\mathcal{M}} f_r\|_{L^2_{\rho_{\mathcal{M}}}}^2 \leq c_{3,\mu} \bigg( 2s + \bigg( \frac{2s}{\lambda} + \frac{1}{s^{d+2}} \bigg) \mathcal{A}(s, \lambda) \bigg).$$

**Proof.** It suffices to notice that both $\lambda \|\vec{f}_\lambda\|_K^2$ and $\mathcal{E}(\vec{f}_\lambda) - 2\sigma_s^2$ are bounded by $\mathcal{A}(s, \lambda)$. □

Next we estimate $\mathcal{A}(s, \lambda)$. We will need the following lemma.

**Lemma B.2.** *Let* $f_r \in C^2(\mathcal{M})$. *There exists* $c_4 > 0$ *such that for* $\vec{f} \in \mathcal{H}_K^p$,

$$\mathcal{E}(\vec{f}) - 2\sigma_s^2 \leq c_4 (s^{d+4} + s^{d+4} \|\vec{f}\|_K^2 + \|\mathscr{D}^* \vec{f} - \nabla_{\mathcal{M}} f_r\|_{L^2_{\rho_{\mathcal{M}}}}^2).$$

**Proof.** Since $\mathcal{M}$ is a submanifold with intrinsic distance, there exists $\delta_0 > 0$ such that $\|\xi - x\| \leq \delta_0$ implies $d_\mathcal{M}(x, \xi) \leq \varepsilon_0$. Denote $\Delta = \{\xi \in \mathcal{M} : \|\xi - x\| \leq \delta_0\}$. Then $\Delta \subset \mathcal{B}_{x,\varepsilon_0}$. So

$$\mathcal{E}(\vec{f}) - 2\sigma_s^2$$
$$\leq \int_\mathcal{M} \int_{\mathcal{B}_{x,\varepsilon_0}} e^{-\|x-\xi\|^2/(2s^2)} (f_r(x) - f_r(\xi) + \vec{f}(x) \cdot (\xi - x))^2 \, d\rho_\mathcal{M}(\xi) \, d\rho_\mathcal{M}(x)$$
$$+ \int_\mathcal{M} \int_{\mathcal{M}\backslash\Delta} e^{-\|x-\xi\|^2/(2s^2)} (f_r(x) - f_r(\xi) + \vec{f}(x) \cdot (\xi - x))^2 \, d\rho_\mathcal{M}(\xi) \, d\rho_\mathcal{M}(x)$$
$$:= J_1 + J_2.$$

It is easy to notice that

$$J_2 \leq \tilde{c}_8 e^{-\delta_0^2/(2s^2)}(1 + \|\vec{f}\|_K^2)$$

for some $\tilde{c}_8 > 0$.

Note that for every $x \in \mathcal{M}$, $\exp_x^{-1}(\mathcal{B}_{x,\varepsilon_0}) \subset B_x(\varepsilon)$. Write $\xi$ in $x$-normal coordinates. Then on $\mathcal{B}_{x,\varepsilon_0}$ there holds $\frac{1}{2}|v|^2 \leq \|x - \xi\|$ and $\sqrt{\det(G)}(v) \leq \frac{3}{2}$. Let $t_1(v)$, $t_2(v)$ and $\tilde{c}_5$ be the same as in the proof of Lemma B.1. We have

$$J_1 \leq 2\tilde{c}_5 \int_\mathcal{M} \int_{\exp_x^{-1}(\mathcal{B}_{x,\varepsilon_0})} e^{-\|x-\exp_x(v)\|^2/(2s^2)} (|t_1(v)|^2 + |t_2(v)|^2) \sqrt{\det(G)}(v) \, dv \, d\rho_\mathcal{M}(x)$$
$$\leq 2\tilde{c}_5 \int_\mathcal{M} \int_{|v|\leq\varepsilon_0} e^{-|v|^2/(4s^2)} (|t_1(v)|^2 + |t_2(v)|^2) \sqrt{\det(G)}(v) \, dv \, d\rho_\mathcal{M}(x)$$
$$\leq \tilde{c}_9((1 + \|\vec{f}\|_K^2) s^{d+4} + s^{d+2} \|\mathscr{D}^*\vec{f} - \nabla_M f_r\|_{L^2_{\rho_\mathcal{M}}}^2)$$

for some $\tilde{c}_9 > 0$.

Combining the estimates for $J_1$ and $J_2$, we finish the proof. $\square$

Lemma B.2 implies the following corollary which bounds $\mathcal{A}(s, \lambda)$ by the functional $\mathscr{K}$.

**Corollary B.2.** *If $f_r \in C^2(\mathcal{M})$, then*

$$\mathcal{A}(s, \lambda) \leq c_5\left(s^{d+4} + s^{d+2} \mathscr{K}\left(\frac{\lambda}{s^{d+2}} + s^2\right)\right)$$

*with $c_5 = \max(c_4, 1)$.*

**Proof.** By Lemma B.2, for every $\vec{f} \in \mathcal{H}_K^p$, there holds

$$\mathcal{E}(\vec{f}) - 2\sigma_s^2 + \lambda\|\vec{f}\|_K^2 \leq c_5(s^{d+4} + s^{d+2}\|\mathscr{D}^*\vec{f} - \nabla_\mathcal{M} f_r\|_{L^2_{\rho_\mathcal{M}}}^2 + (\lambda + s^{d+4})\|\vec{f}\|_K^2).$$

The conclusion follows by taking infimum over $\vec{f} \in \mathcal{H}_K^p$. $\square$

One easily sees that Propositions B.2 and B.3 follow from combining Corollaries B.1 and B.2.

***Remark B.1.*** We remark that the approximation error estimate in Proposition B.2 converges to 0 as $s \to 0$ if $\mathscr{K}(t) = \mathrm{O}(t^\beta)$ with $\beta > \frac{1}{2}$. The result may be improved by using functional analysis techniques; see, for example, [16]. However, it seems those techniques require the explicit functional expression of $\vec{f}_\lambda$, which is available only in the regression case. The proof method we provide here is not as powerful for the regression case but it is more general and can be applied even in cases where $\vec{f}_\lambda$ only exists implicitly, such as the classification setting.

Now we prove Theorems 3.1 and 3.2.

**Proof of Theorem 3.1.** First note that $\mathscr{K}(t)$ is increasing with $t$. Since $s \leq 1$ and $\lambda = s^{d+2+2\theta} \leq 1$, we have $\mathscr{K}(\frac{\lambda}{s^{d+2}} + s^2) = \mathscr{K}(s^{2\theta} + s^2) \leq \mathscr{K}(2s^{2\theta})$. Then by Corollary B.2,

$$\lambda \|\vec{f}_\lambda\|_K^2 \leq \mathcal{A}(s, \lambda) \leq c_5(s^{d+4} + s^{d+2}\mathscr{K}(2s^{2\theta})).$$

Plugging into the sample error estimate in Proposition B.1 gives

$$\|\vec{f}_D - \vec{f}_\lambda\|_K \leq \frac{C' \log(2/\delta)}{\sqrt{m\lambda}}(1 + s^{-\theta}\sqrt{\mathscr{K}(2s^{2\theta})})$$

with $C' = C_1 + (C_2 + \frac{1}{\log 2})\sqrt{c_5}$. By Proposition B.2

$$\|\vec{f}_\lambda - \nabla_\mathcal{M} f_r\|_{L^2_{\rho_\mathcal{M}}}^2 \leq 3C_3(s^\theta + s^{-\theta}\mathscr{K}(2s^{2\theta})).$$

Combining these two estimates, we draw the conclusion with the constant $C_\rho = \max\{(C')^2, 3C_3\}$. $\square$

Note that if $\mathcal{M}$ is a compact domain in $\mathbb{R}^p$, then $d = p$, $\mathscr{D} = \mathscr{D}^* = I$, and $\nabla_\mathcal{M}$ is the usual gradient operator. In this case $\mathscr{K}(t) = \mathrm{O}(t)$ if $\nabla f_r \in \mathcal{H}_K^p$. The rate in Corollary 3.1 becomes $\mathrm{O}(n^{-\theta/(2p+4+5\theta)})$, which is of the same order as that derived in [16]. This implies that our result reduces to the Euclidean setting when the manifold is a compact domain in the Euclidean space.

**Proof of Theorem 3.2.** By $s \leq 1$ and $\lambda = s^{d+3}$ we obtain $\mathscr{K}(\frac{\lambda}{s^{d+2}} + s^2) = \mathscr{K}(s + s^2) \leq \mathscr{K}(2s)$. Then by Corollary B.2,

$$\lambda \|\vec{f}_\lambda\|_K^2 \leq \mathcal{A}(s, \lambda) \leq c_5(s^{d+4} + s^{d+2}\mathscr{K}(2s)).$$

Plugging into the sample error estimate in Proposition B.1 gives

$$\|\vec{f}_D - \vec{f}_\lambda\|_K \leq \frac{C' \log(2/\delta)}{\sqrt{m}\lambda}(1 + \sqrt{s^{-1}\mathscr{K}(2s)})$$

with $C' = C_1 + (C_2 + \frac{1}{\log 2})\sqrt{c_5}$. By Proposition B.3,

$$\|\vec{f}_\lambda - \nabla_\mathcal{M} f_r\|^2_{L^2_{\rho_\mathcal{M}}} \leq 3C_{3,\mu}(s + \mathscr{K}(2s^2)).$$

The conclusion follows by combining the above two estimates. □

# Acknowledgements

The work was partially supported by the Research Grants Council of Hong Kong [Project No. CityU 104007] and a Joint Research Fund for Hong Kong and Macau Young Scholars from the National Science Fund for Distinguished Young Scholars, China [Project No. 10529101]. The work was also partially supported by NSF Grant DMS-07-32260, NIH Grant R01 CA123175-01A1 and NIH Grant P50 GM 081883.